\documentclass[12pt]{article}
\usepackage{amsxtra,amssymb,amsthm,amsmath,latexsym}

\theoremstyle{plain}

\numberwithin{equation}{section}

\def\qed{{\hfill $\Box$}}

\def\ve{{\varepsilon}}

\def\oH{{\overset{\circ}{H}}}
\def\oH1{{\overset{\circ}{H}\kern-.02in{}^1}}

\def\bee{\begin{equation*}}
\def\eee{\end{equation*}}
\def\be{\begin{equation}}
\def\ee{\end{equation}}

\begin{document}
\title{   Discrete L'Hospital's rule 
}

\author{
A.G. Ramm\\
 Mathematics Department, Kansas State University, \\
 Manhattan, KS 66506-2602, USA\\
ramm@math.ksu.edu\\}

\date{}

\maketitle\thispagestyle{empty}

\footnote{2000 Math subject classification:
26A24, 26D15 }


\section{Introduction}\label{S:1}

The aim of this paper is to formulate discrete analog of L'Hospital's
rule and describe some of its applications.
While the usual L'Hospital's rule is taught to all undergraduate students
studying calculus, its discrete analog apparently was not in the 
literature. Since the  L'Hospital's rule proved to be very useful in many 
applications, one may think that its  discrete analog will also be useful.
We start by stating the usual, well-known L'Hospital's rule, so that 
the reader could see the similarities in the formulation of this rule and 
its discrete analog. Then we formulate this  discrete analog.
In Section 2 we prove the  discrete analog and illustrate it by some 
examples. 

After this note has been written, the author learned that in \cite{F}, 
p.67, there is a result of O.Stolz, which is the same as Theorem 2 below.
Our proof is slightly different from the proof in \cite{F}.
The application, given in our Example 2 is of interest in the 
theory of the dynamical systems method (\cite{R}).

A version of the usual  L'Hospital's rule is the following theorem.

{\bf Theorem 1.} {\it Assume that:

1) the functions $F$ and $G$ are continuously differentiable on the 
interval $I=(a,a+h)$, where $h>0$ and $a$ are real numbers, 
$f:=F'$, $g:=G'$, and 
\be\label{e1.1}
\lim F=\lim G=\infty, \quad \lim:=\lim_{x\to a, x>a},
\ee 
2)  
\be\label{e1.2}
\lim \frac f g=L, \quad g(x)\neq 0 \quad \forall x\in I.
\ee
Then there exists the limit:
\be\label{e1.3}
\lim \frac F G=L.
\ee
}

The proof of Theorem 1 can be found in any calculus text and does not need 
a reference.

Let us now formulate the discrete analog of the above theorem.

{\bf Theorem 2.} {\it Let $f_j>0$ and $g_j>0$  be sequences of numbers,
$F_n:=\sum_{j=1}^n f_j$, $G_n:=\sum_{j=1}^n g_j$. 

Assume that:
\be\label{e1.4}
\lim_{n\to \infty} F_n=\lim_{n\to \infty} G_n =\infty,
\ee
and 
\be\label{e1.5}
\lim_{n\to \infty} \frac {f_n}{g_n}=L.
\ee
Then 
\be\label{e1.6}
\lim_{n\to \infty} \frac {F_n}{G_n}=L.
\ee
}
The similarity of Theorems 1 and 2 is obvious.

{\bf Remark:} {\it One can write equation (1.5) as
\be\label{e1.7}
\lim_{n\to \infty} \frac {F_n-F_{n-1}}{G_n-G_{n-1}}=L.
\ee 
Thus, the role of the derivative of $F$ is played by the 
difference $F_n-F_{n-1}$.
}

\section{Proofs}\label{S:2}
Fix an arbitrary small $\ve>0$. Using 
assumption (1.5), find $M:=M(\ve)$, such that
\be\label{e2.1}
 L-\ve<\frac {f_j}{g_j}<L+\ve,  \quad
\forall j>M. \ee 
Denote $F_{nM}:=\sum_M^n f_j$, and define $G_{nM}$  similarly.
Using assumption (1.4), find $N:=N(\ve)$, such that
\be\label{e2.2}
 \frac {F_M}{F_{nM}}<\ve,   \quad \frac 
{G_M}{G_{nM}}<\ve,\quad
\forall n>N. \ee
Now one gets:
\be\label{e2.3}
 \frac {F_n}{G_n}= \frac {F_M+F_{nM}}{G_M+G_{nM}}=\frac{F_{nM}}{G_{nM}}
\frac  {1+\frac {F_M}{F_{nM}}}{1+\frac {G_M}{G_{nM}}},
\ee
and 
\be\label{e2.4}
 1+\ve_1:=\frac {1-\ve}{1+\ve}<
\frac  {1+\frac {F_M}{F_{nM}}}{1+\frac {G_M}{G_{nM}}}<\frac 
{1+\ve}{1-\ve}:=1+\ve_2,
\ee
where $\ve_1=O(\ve)$ and  $\ve_2=O(\ve)$, as $\ve\to 0$.
Using assumption (1.5), one gets
\be\label{e2.5}
 L-\ve\leq \min_{j\geq M} \frac {f_j}{g_j}\leq
\frac {F_{nM}}{G_{nM}}\leq \max_{j\geq M} 
\frac{f_j}{g_j}\leq  L+\ve.
\ee
Since $\ve>0$ is arbitrarily small, equation (1.6) follows
from relations (2.2)-(2.5).
Theorem 2 is proved. \qed

Consider examples of applications of Theorem 2.

{\it Example 1.} By Theorem 2, one has $$ \lim_{n\to \infty}
\frac {\sum_{j=1}^n \frac {j^m}{1+j^{m+1}}} {\sum_{j=1}^n
\frac {j^p}{1+j^{p+1}}}=\lim_{n\to \infty}
\frac{{n^m}/{(1+n^{m+1})}}{{n^p}/{(1+n^{p+1}})}=1. $$ 

{\it Example 2.} In applications the following differential
inequality is used (see, e.g., \cite{R}): 
\be\label{e2.6} g_{n+1}\leq
(1-a_n)g_n+b_n, \quad n=1,2,3,..... 
\ee 
Assume that
\be\label{e2.7} 0<a_n<1, \quad\lim_{n\to \infty} \frac
{b_{n-1}}{a_n}=0,\quad \sum_{n=1}^\infty a_n=\infty.
\ee 
Note that assumptions (2.7) imply $\lim_{n\to \infty}b_n=0$.

Using assumptions (2.7) one can apply Theorem 2 and conclude 
that 
\be\label{e2.8} 
\lim_{n\to \infty}
\sum_{k=1}^{n-1} b_k\prod_{j=k+1}^n(1-a_j)=0. \ee 

This result
implies, that  $\lim_{n\to \infty}g_n=0$ under the 
assumptions (2.7), where $g_n$ is a sequence solving 
inequality (2.6).

Let us discuss in detail the application of Theorem 2 in
this example.

 From (2.6) by induction one gets:
\be\label{e2.9}
g_{n+1}\leq b_n+\sum_{k=1}^{n-1}b_k \prod_{j=k+1}^n
(1-a_j)+g_1\prod_{j=1}^n (1-a_j).
\ee
Assumption (2.7) implies that 
\be\label{e2.10} \lim_{n\to \infty}b_n=0\quad  \text {  and 
} 
\lim_{n\to 
\infty}g_1\prod_{j=1}^n (1-a_j)=0.
\ee
Let us write the term $J_n:=\sum_{k=1}^{n-1}b_k 
\prod_{j=k+1}^n(1-a_j)$ in the form:
$J_n= \frac {\sum_{k=1}^{n-1}b_k\prod_{j=1}^k(1-a_j)^{-1}}
{\prod_{j=1}^n (1-a_j)^{-1}} $.
We want to apply Theorem 2 in order to prove that
 \be\label{e2.11}\lim_{n\to \infty}J_n=0.\ee 
The denominator in $J_n$ tends to infinity. If the numerator
in  $J_n$ is bounded, then (2.11) follows. If
this numerator tends to infinity, then one has 
assumption (1.4) satisfied. To check assumption (1.5) with 
$L=0$, one 
calculates the limit:
$$
\lim_{n\to \infty} \frac {b_{n-1}\prod_{j=1}^{n-1} (1-a_j)^{-1}}
{\prod_{j=1}^n (1-a_j)^{-1}[1-(1-a_n)]}=lim_{n\to \infty} 
\frac {b_{n-1}(1-a_n)}{a_n}=0.$$
At the last step asssumption (2.7) was used.
So, Theorem 2 yields the desired conclusion (2.11).
The discussion of Example 2 is completed.\qed


\begin{thebibliography}{1000} 

\bibitem{F} Fikhtengolts, G., Course of differential and integral 
calculus, vol.1, Fizmatgiz, Moscow, 1962.

\bibitem{R} Ramm, A.~G.~,  
 Dynamical systems method for ill-posed equations with
monotone operators, Comm. in Nonlinear Sci. and Numer.
Simulation, 10, N2, (2005).








\end{thebibliography}
\end{document}